\documentclass[a4paper,12pt]{amsart}
\usepackage{tabularx}
\usepackage[dvipdfm]{graphicx}
\usepackage{amsmath}
\usepackage{amssymb}
\usepackage{amsbsy}
\usepackage{amsgen}
\usepackage{amsopn}
\usepackage{amscd}
\usepackage{amstext}
\usepackage{amsthm}

\theoremstyle{definition}
\newtheorem{dfn}{Definition}[section]

\theoremstyle{plain}
\newtheorem{thm}{Theorem}[section]
\newtheorem{lem}{Lemma}[section]
\newtheorem{cor}{Corollary}[section]
\newtheorem{prop}{Proposition}[section]

\theoremstyle{remark}
\newtheorem{rem}{Remark}[section]

\theoremstyle{definition}
\newtheorem{ex}{Example}[section]
\newtheorem{acknow}{Acknowledgment}

\newtheorem{conv}{Conventions}

\begin{document}

\title{Invariants via word for curves and fronts}
\author{Noboru Ito}

\maketitle

\begin{abstract}  
We construct the infinite sequence of invariants for curves in surfaces by using word theory that V. Turaev introduced.  
For plane closed curves, we add some extra terms, e.g. the rotation number.  From these modified invariants, we get the Arnold's basic invariants and some other invariants.  We also express how these invariants classify plane closed curves.  In addition, we consider other classes of plane curves: long curves and fronts.  
\end{abstract}

\section{Introduction}\label{jyo}

The object of the paper will be to construct 
some invariants of plane curves and fronts,  and so it is to show one of the method for applying word theory to plane curves and fronts.  
V. Turaev introduces word theory  (\cite{turaev1}, \cite{turaev2}, \cite{turaev3}).  
We can consider that this word theory is effective in two view points as follows.  
\begin{enumerate}
\item Word is the universal object of  knot, curve, etc.  \label{item:n1}
\item We can treat knots and curves collectively and
 algebraically, so that we can systematically study in these invariants themselves and relationships among them.  \label{item:n2}
\end{enumerate}
In terms of (\ref{item:n1}), Turaev applies  topological methods (Reidemeister move, 
homotopy, etc.) to a semigroup consisting
 of letters, so that creates word which has property as (\ref{item:n1}) \cite{turaev1}.  
In terms of (\ref{item:n2}), Turaev considers
equivalent classes of words corresponding to knots  or curves, 
and constructs invariants of knots, for 
example, Jones polynomial and 
$\underline{\alpha}$-kei which is similar to kei
 for knots \cite{turaev2}.  
 
For immersed plane closed curves, H. Whitney classified plane closed curves regular homotopically 
by winding number, which is also called index or rotation number \cite{W}.  
Long afterword, V. I. Arnold created three
 basic invariants of plane closed curves $J^+, J^-, St$ by the similar method to knots of V. A. Vassiliev \cite{V} and classified plane closed curves which have same index  (\cite{arnold1}, \cite{arnold2}).  
Arnold also obtained a natural generalization of $J^+$ to fronts (\cite{arnold1}, \cite{arnold3}).  
Relating to this, M. Polyak systematically reconstructed the Arnold's  basic invariants via Gauss diagram and related basic invariants to the Vassiliev invariant \cite{polyak1}.      

In this paper, by using word theory, 
we will reconstruct the Arnold's basic invariants and construct 
some other invariants for plane closed curves,
 long curves, and fronts.  
We also express how these invariants classify these plane curves and fronts.  

The outline of each section is as follows.  
In Section \ref{id}, we will compose
 invariants $\{I_n\}$ (`\textbf{i}'nvariant of degree `\textbf{n}') of curves on a surface.  
In Section \ref{plane}, we will construct invariants of plane closed curves $CI_n$
 (`\textbf{c}'losed curve `\textbf{i}'nvariant of degree`\textbf{n}')
for $I_n$.  
$CI_2$ has the same strength as the Arnold's basic invariants.  $CI_3$  is independent of $CI_2$.  
There is an example that two curves take the
 same values of
 index, the Arnold's basic invariants and HOMFLY polynomial of immersed plane closed curves \cite{CGM} but take different values of $CI_3$.  
In Section \ref{long}, Section \ref{front}, we study in long curves and fronts by using
 the similar technique.  
    
\begin{conv}
In this paper, all surfaces and curves are
 oriented.  
For a given surface $\varphi$, 
a {\itshape closed curve} (resp. {\itshape long curve}) is an immersion
 : $S^1$ (resp. $\mathbf{R}$) $\to$ $\varphi$ (resp. $\mathbf{R}^2$) where 
all of the singular points are transversal 
double points.  
A {\itshape front} is  an immersion
 : $S^1$ $\to$ $\mathbf{R}^2$ 
with the coorientation (defined in Sect. \ref{dfn-front}) where 
all of the singular points are transversal 
double points or cusps.  
(We will precisely define 
a {\itshape front} in Sect. \ref{dfn-front}
.  ) 
A {\itshape curve} is a closed curve, 
a long curve, or a front.  
A {\itshape smooth curve} is a closed curve 
or 
a long curve.  
When a curve stands for a closed curve
 or a front, 
a {\itshape base point} is a point on
 the curve except on the double points
 and the cusps.  
A {\itshape pointed curve} is a closed curve or a front endowed with a base point.  
Winding number (rotation number) is
 called {\itshape index} in this paper.  
\end{conv}

\begin{acknow}
I am grateful to Professor Jun Murakami for 
giving me numerous fruitful comments.  
This paper is nothing but trying to answer in my way to his question hitting the mark : ``Can we apply word theory to plane curve theory?"  
\end{acknow}

\section{Invariants $\{I_n\}$}\label{id}

In this section we equip Turaev's word 
to construct the sequence of invariants for pointed  surface curves.  

\subsection{Turaev's word}\label{word}
We follow the notation and terminology of \cite{turaev2}.  
An {\itshape alphabet} is a set and its elements  are called {\itshape letters}.  
A {\itshape word of length} $m \ge 1$ in an alphabet $\mathcal{A}$ is a mapping 
$\hat{m}=\{1, 2, \ldots, m-1, m\} \to \mathcal{A}$.  
A word $w : \hat{m} \to \mathcal{A}$ is encoded by the sequence of letters 
$w(1)w(2) \cdots w(m)$.  
Two words $w$, $w'$ are {\itshape isomorphic} 
if there is a bijection $w' = f w.  $ 

A word $w$ is called a {\itshape Gauss word} 
if each element of $\mathcal{A}$ is the image of precisely two elements of $\hat{m}$.     
For an alphabet $\alpha$, 
an {\itshape $\alpha$-alphabet} $\mathcal{A}$ is an  alphabet endowed with a 
projection $|~| : \mathcal{A} \to \alpha$.  
An {\itshape \'{e}tale word over $\alpha$} is a pair $(\mathcal{A}, w)$
 where $\mathcal{A}$ is $\alpha$-alphabet and $w:\hat{m} \to \mathcal{A}.  $
In this paper, we only treated \'{e}tale word 
 $(\mathcal{A}, w)$ where $w$ is a surjection.  
In particular, a {
\itshape nanoword over $\alpha$} is an \'{e}tale word $(\mathcal{A}, w)$ over $\alpha$
 where $w$ is a Gauss word.  
For $(\mathcal{A}, w)$, we admit that we use the simple description `$w$' if this $w$ 
means $(\mathcal{A}, w)$ clearly.   
An {\itshape isomorphism} of $\alpha$-alphabets $\mathcal{A}_1$, $\mathcal{A}_2$ is a bijection 
$f : \mathcal{A}_1 \to \mathcal{A}_2$ endowed with $|f(A)|=|A|$ 
for all $A \in \mathcal{A}_1$.  
Two nanowords $(\mathcal{A}_1, w_1)$, $(\mathcal{A}_2, w_2)$ over $\alpha$
 are {\itshape isomorphic} if there is an isomorphism $f$ of
 $\alpha$-alphabets $\mathcal{A}_1$, $\mathcal{A}_2$ such that 
$w_2 = f w_1.  $

Until Sect. \ref{long}, we denote by ${\mathit
{sign}}$ the projection $|~|$.  
We also define an alphabet $\alpha$, an involution $\tau$, and a set $S$ by 
\[
\alpha = \{-1, 1\}, \quad 
\tau : -1\mapsto 1, \quad 
S = \{ (-1, -1, -1),  (1, 1, 1)\}.  
\]until Sect. \ref{long}.  

The following fundamental theorem is established by Turaev \cite{turaev2}.  

\begin{figure}[htbp]
\begin{center}
\includegraphics[width=13cm]{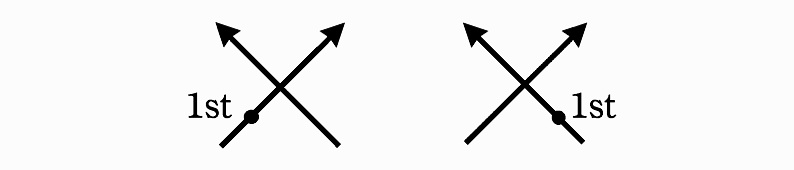}

{\Large $-1$\qquad\qquad\qquad$1$}
\caption{Two patterns of the double points correspond to two letters in $\alpha$.  }
\label{yajirusi}
\end{center}
\end{figure}

\begin{thm}\label{curve} {\rm (Turaev)}
Every pointed closed curve is represented as a nanoword.   
\end{thm}
\begin{proof}
 For a given pointed closed curve which has precisely $m$ double points, 
we name the double points $A_1, A_2, \ldots, A_m$ along the curve orientation from the base point.  
Each point precisely corresponds to either $-1$ or $1$ in Figure \ref{yajirusi}.  
\end{proof}

\begin{rem}
Theorem\ref{curve} implies that there is a 
mapping from a smooth curve $\Gamma$ on a surface
 to a nanoword.  
 Conversely, there is a mapping from a nanoword to a smooth curve on a surface.  
In other words, we determine a unique smooth 
curve $\Gamma$ and a unique surface 
$\Sigma$ on which $\Gamma$ is by using the 
following Theorem\ref{word-curve}.  
\end{rem}

\begin{thm}\label{word-curve}
\rm{(Turaev)}
Let $w$ be a nanoword of length $n$ and
 $g(\Gamma)$ the minimum genus of the compact
 surface $\Sigma$
 without boundary such that a pointed smmoth
 curve $\Gamma$ is on $\Sigma$.  
There is a mapping $w \mapsto \Gamma$ and 
\[g(\Gamma) = 1 + \frac{n - card(\bar{n}/t\theta)}{2}.  \]
\end{thm}

The construction of the mapping $w \mapsto \Gamma$ is well known and therefore 
this is omitted (cf. \cite{turaev5}, \cite{turaev4}, \cite{viro}).  
We can calculate $g(\Gamma)$ by the proof of
Theorem 9.1.1 
in \cite{turaev4} for Turaev's {\itshape chart}.

\subsection{Construction of invariants $I_n$}
\label{constIn}

For every Gauss word $v$
 and every nanoword $w$, 
we determine a number 
$\left\langle v, w\right\rangle$.  
When a nanoword $(\mathcal{A}_w, w)$ over 
$\alpha$
 is given, we consider a sub-word $v'$ of $w$.  
If a sub-word $v'$ is Gauss word, 
we can naturally consider  the nanoword $ (\mathcal{A}_{v'}, {v'}) $ over $\alpha$ such that $\mathcal{A}_{v'} \subset \mathcal{A}_w$.  
Therefore for every nanoword $(\mathcal{A}_w, w)$ over
 $\alpha$ and for every Gauss word $v$, we can define 
the mapping by
\[\left\langle v, w\right\rangle :=\sum_{\text{a sub-word $v'$ of $w$ isomorphic to $v$}} \quad
{\prod}_{A \in A_{v'}} \mathit{sign}A.  \]
Let $W_n$ be the free $\mathbf{Q}$-module generated by the set of 
all of the isomorphism class of the Gauss words where each length of the Gauss word is $2n$.  
For a given integer $d$, let $N_d$ be the free $\mathbf{Q}$-module generated by the set of nanowords over $\alpha$ where each length 
of the nanoword is less than $2d+1$.  
Expanding $\left\langle v, w\right\rangle$  bilinearly, 
we can make a bilinear mapping $\left\langle, \right\rangle$ 
from $W_n \times N_d$ to $\mathbf{Q}.  $  
For an arbitrary surface, 
let $w_\Gamma$ stand for a word which is determined
by a curve $\Gamma$ on the surface.  

\begin{thm}\label{I_d}
The following $\{I_n\}$ (invariant of degree $n$) is the sequence of  surface isotopy invariants for
 pointed curves on a surface.  

\[I_n (\Gamma)=\left\langle \sum_k x_k v_k, w_\Gamma\right\rangle
\qquad  (n\in \mathbf{N})\]

where $\{v_k\}$ is the basis of $W_n$ and 
each $x_k$ is a parameter.  
\end{thm}

\begin{proof}
By using Theorem \ref{curve}, the way  of constructing 
$\langle , \rangle$ implies this theorem.  
\end{proof}

\subsection{Generalization of $I_n$}\label{Gin}

We can generalize $I_n$ by introducing  {\itshape a dimension of a letter}.  

\begin{dfn} (a dimension of a letter, a dimension of a word)\label{dim}
Let $X$, $\dot{X}$ be letters. For an arbitrary letter $A$, we denoted by
 $d (A) \in \{1, 2\}$
 {\itshape a dimension of a letter $A$}
 which is defined by the following.  
Let $X$ be a $1$-dimensional letter where
 $d(X)=1$ and let $\dot{X}$
 a 2-dimensionl letter where
 $d(\dot{X})=2$.  
Next, let $\mathcal{A}$ be an alphabet.  
For every word
 $w : \hat{m} \to \mathcal{A}$, 
the dimension $d(w)$ of word $w$
 is defined by $d(w) := \sum_{X \in \mathcal{A}} d(X)$.  
\end{dfn}
 
The concept of the dimension of word affect on the module `$W_n$'.  That is why we must redefine `$W_n$'.  

The following {\itshape word space} $\mathcal{W}_n$ is the canonical generalization of 
`$W_n$' defined in Sect. \ref{constIn}.  

\begin{dfn} (word space)\label{gokukan}
The {\itshape word space} $\mathcal{W}_n$ {\itshape of degree $n$} is the free $\mathbf{Q}$-module generated by the set of  all of the $n$-dimensional Gauss words
 which may contain $2$-dimensional letters.  
\end{dfn}

Replacing $W_n$ defined in Sect. \ref{constIn} with $\mathcal{W}_n$, we can easily check that the similar results are established and can easily generalize Sect. \ref{constIn}.     

For an arbitrary $(\mathcal{A}_w, w)$ and $v
 \in \mathcal{W}_n$, 
we think of $(\mathcal{A}_{v'}, {v'})$ 
where $v'$ is a sub-word of $w$ which is  isomorphic to $v$ and $\mathcal{A}_{v'} \subset
{A}_{w}$.  
We {\itshape redefine} $\left\langle,  \right\rangle$
 by 
\[\left\langle v, w\right\rangle :=\sum_{\text{a sub-word $v'$ of $w$ isomorphic to $v$}} \quad \prod_{A \in A_v} \left(\mathit{sign}A\right)^{d(A)}.  \]
 
\begin{cor}\label{gin}
The following $\{GI_n\}$ (generalized $I_n$) is the sequence of  surface isotopy invariants for
 pointed curves on a surface.  

\[GI_n (\Gamma)=\left\langle \sum_k x_k v_k, w_\Gamma\right\rangle
\qquad  (n\in \mathbf{N})\]

where $\{v_k\}$ is the basis of the word space $W_n$ and 
each $x_k$ is a parameter.  
\end{cor}

\begin{rem}
For every plane closed curve $\Gamma$, 
$\left \langle \dot{X}\dot{X} ,
 w_\Gamma \right \rangle $ is the number of the double points for $\Gamma$ \cite{polyak1}.  
\end{rem}

\section{Application of $I_n$ to plane closed curves}\label{plane}

We will consider an application of $I_n$ to plane closed curve.  
We can apply word theory to plane curve theory because word is universal for knots  and curves.  
In this section we will apply it to plane closed curve for example.  
Plane curves are not only fundamental objects but also proper objects to think of some various applications of word theory.  
In fact, we can apply word theory to closed curves, long curves, and fronts (Sect. \ref{long}, \ref{front}).  
When we apply word theory to plane closed curves, we get some invariants $CI_n$ (`\textbf{c}'losed curve `\textbf{i}'nvariant degree $n$).  
In order to construct $CI_n$, to add to Turaev's word, we need one more material 
about plane curve theory : 
the Arnold's basic invariants defined 
in the next subsection.  

\subsection{The Arnold's basic invariants}
\label{arnold-invariant}  
We consider regular homotopy classes of plane curves.  
Let us rewrite the Arnold's the invariants via Turaev's word theory.     
To redefine the Arnold's basic invariants \cite{arnold2}, 
we define {\itshape elementary moves} that are local moves (Figure \ref{2move}, \ref{3move})
 of plane 
curves apart from a base point.  

\begin{dfn} (elementary move) 
Let $x, y, z$ be words that consist of the letter in $\mathcal{A}$ where $\mathcal{A} \cap \{A, B\}= \emptyset$.  
Elementary move $I\!I^+$ and elementary move $I\!I^-$ (Figure \ref{2move}) are defined by 
\[I\!I^+ : (\mathcal{A}, xyz)\to  (\mathcal{A}\cup\{A, B\}, xAByABz)  ~\text{if}~ \tau ({\mathit sign}A)= {\mathit sign}B, \]
\[I\!I^- : (\mathcal{A}, xyz)\to  (\mathcal{A}\cup\{A, B\}, xAByBAz) ~\text{if}~ \tau ({\mathit sign}A)= {\mathit sign}B.  \]

Let $x, y, z, t$ be words that consist of the letter in $\mathcal{A}-\{A,B,C\}$.   
Elementary move $I\!I\!I$ (Figure \ref{3move}) is defined by 
\[(\mathcal{A} , xAByACzBCt) \to 
 (\mathcal{A}, xBAyCAzCBt)~\text{for}~({\mathit sign}A, {\mathit sign}B, {\mathit sign}C)\in S.  \] 
The positive elementary moves is the above direction, the negative elementary move
 is the inverse direction.  
\begin{figure}[htbp]
\begin{center}
\includegraphics[width=10cm]
{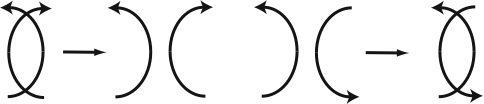}
\caption{Positive elementary move $I\!I^+$ (left figure) and Positive elementary move $I\!I^-$ (right figure).  }
\label{2move}
\end{center}
\end{figure}

\begin{figure}[htbp]
\begin{center}
\includegraphics[width=20cm,height=3cm,keepaspectratio]
{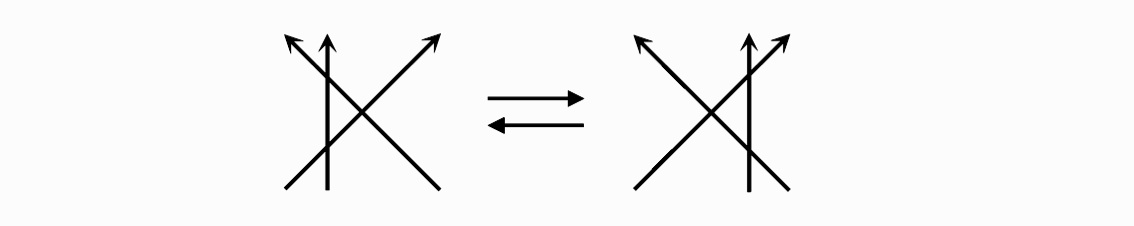}
\caption{ Elementary move $I\!I\!I$.  If the order of three branch are 1, 2, 3 from the left, 
the direction of the positive move is from the left figure to the right figure.  
If the order of three branch are 3, 2, 1 from the left, the direction of the positive move is from the right figure to the left  figure.  }
\label{3move}
\end{center}
\end{figure}
\end{dfn}

For this elementary moves, 
the Arnold's basic invariants $J^+, J^-, St$ are 
invariants of curve can be defined by following  (cf. \cite{arnold2}).  

\begin{dfn} (the Arnold's basic invariants)
\label{basicinvariant}
$J^+$ is increased by 2 under positive elementary move $I\!I^+$ but not change under the other, 
$J^-$ is decrease by 2 under positive elementary move $I\!I^-$ but not change under the other, 
$St$ is increased by 1 under positive elementary move $I\!I\!I$ but not change under the other, 
and 
satisfy the following conditions 
\[J^+ (K_0)=0, \quad J^- (K_0)=-1, \quad St (K_0)=0, \]
\[J^+ (K_{i+1})=-2i, \quad J^- (K_{i+1})=-3i, \quad St (K_{i+1})=i \]
for the base curves $\{K_i\}$ defined by Figure 4.  
\begin{figure}[htbp]
\begin{center}
\includegraphics[width=15cm]
{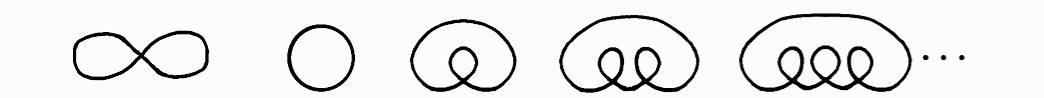}
\begin{flushleft}
\qquad \qquad \ $K_0$
\qquad \qquad \ \  $K_1$
\qquad \quad \   $K_2$
\qquad \qquad \  $K_3$
\qquad \qquad \quad \ $K_4$
\caption{Base curve $\{K_i\}$.  }
\end{flushleft}
\label{base-curve}
\end{center}
\end{figure}
\end{dfn}

\subsection{Construction of invariants $\overline{CI}_n$}\label{const}

In this subsection, we compose a mapping 
$\left[, \right]$ to construct invariants
 $\overline{CI}_n$ of plane closed curves.  
To add to $\left\langle, \right\rangle$, 
we equip {\itshape cyclic equivalent} 
to construct a mapping 
$\left[, \right].  $

\begin{dfn} (cyclic equivalent)\label{jyunkaidouti}
Let $x$ be $w (2)$ $\cdots$ $w (m)$ of $w$
 $=$ $w (1)$$w (2)$ $\cdots$ $w (m)$, 
for two arbitrary Gauss words $w, w' \in W_n$, the  relation $\sim$ is defined by
\[
w \sim w' \stackrel{\mathrm{def}}{\Longleftrightarrow} 
w=Ax~\text{and}~w'=-xA.  \]  
This relation $\sim$ is called {\itshape cyclic equivalent}.  
\end{dfn}

The cyclic equivalent is equivalent relation.    
Let $\overline{W}_n$ be a module consisting of cyclic equivalent classes
 of the elements of $W_n$ (defined by Sect. \ref{constIn}).  
For $[w] \in \overline{W}_n$, the number of the residue system of $[w]$
 is even.  That is because $w \sim -w$ implies $w \sim 0$ if this number is odd.  

The mapping \[\left[, \right] : \overline{W}_n \times N_d \to \mathbf{Q}\]
is defined by
\[\left[[v], w \right]:=
\left<v_1+v_2+ \cdots +v_{2l} , w \right> \] 
where $v_1, v_2, \ldots, v_{2l}$ consist of all elements of $[v]$. 

\begin{prop}\label{[]}
Let $\Gamma$ be an arbitrary curve.  
For every $[w] \in \overline{W}_n$, 
$[[w], w_\Gamma]$ is a surface isotopy invariant of curves.  
\end{prop}

\begin{proof}
Base point move 
 (Figure \ref{basepoint-move}) is , that is to say,  
 to replace $(\mathcal{A}, w_\Gamma)=AxAy 
$ with $xA_\tau yA_\tau, {\mathit sign}{A_\tau} =\tau ({\mathit sign}A)$ 
where $x, y$ are consist of the letters of $\mathcal{A}-\{A\}$.    
\begin{figure}[htbp]
\begin{center}
\includegraphics[width=15cm,height=2cm,keepaspectratio]
{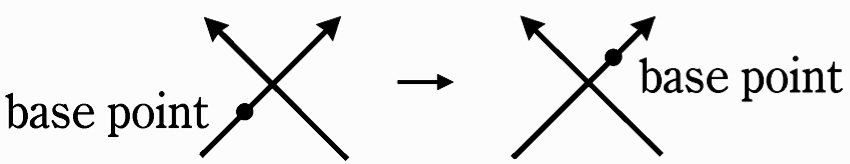}
\caption{Base point move.  }
\label{basepoint-move}
\end{center}
\end{figure}
Therefore under the base point move, a part of 
$\left<AxAy, w_\Gamma \right>$ multiplied -1 
is added to 
$\left<xAyA, w_\Gamma \right>$.  
By definition of cyclic equivalence, and 
by the new numbering

$v_1, v_2, \ldots , v_{2l}$ if necessary, we can have 
\[v_1 \stackrel{\text{base point move}}{\to}
v_2\stackrel{\text{base point move}}{\to} 
\cdots \stackrel{\text{base point move}}{\to}
v_{2l} \stackrel{\text{base point move}}{\to} v_1.  \]  
therefore, the value of 
$\left<v_1+v_2+ \cdots +v_{2l} , w_\Gamma \right>
$ is not change by base point move.  
\end{proof}

\begin{cor}
The following $\{\overline{I}_n\}$ is 
the sequence of surface isotopy invariants
for 
closed curves on a surface.  

\[\overline{I}_n(\Gamma) := \left[ \sum_k x_k [v_k], w_\Gamma\right]
\qquad  (n\in \mathbf{N}).  \]
$\{[v_k]\} $ is the base of $\overline{W}_n$, and  
each $x_k$ is a parameter.  

In particular, 
for every closed plane curve, 
the following $\{\overline{CI}_n\}$ is the sequence of plane isotopy invariants.  
\[\overline{CI}_n (\Gamma) := \overline{I}_n(\Gamma) + f(i) \qquad  (n\in \mathbf{N})\]
where f is function of index i.  
\end{cor}

Next subsection, we introduce $CI_2$ and $CI_3$
which are made of $\overline{CI}_2$ and $\overline{CI}_3$.

\subsection{$CI_2$ and $CI_3$}

For every curve $\Gamma$, let $i$ be index and 
$n$ the number of the double points, 
we define $CI_2$ by 
\[CI_2 (\Gamma; s, t, u) := s n + \left \langle t XXYY - t XYYX + u XYXY, w_\Gamma\right\rangle
+\frac{t}{2}-\frac{t}{2}i^2.  \]  

\begin{rem}
$\overline{CI}_2 (\Gamma)
 = \left[t[XXYY], w_\Gamma \right]$, and
  then, we have 
\[CI_2 (\Gamma;  s, t, u)
 = \overline{CI}_2 (\Gamma) +
 s n + \left \langle u XYXY, w_\Gamma\right\rangle
+\frac{t}{2}-\frac{t}{2}i^2.  \]
\end{rem}

M. Polyak proved that 
$\left\langle
XYXY, w_\Gamma\right\rangle$ does not depend on the choice of a base point 
(cf. Theorem 1 proof in \cite{polyak1}).    
Therefore 
the invariant $CI_2$ is well-defined.  
$CI_2$ is also not depend on the orientation of the curve 
$\Gamma$ because this formula is symmetric.  
$CI_2 (\Gamma; s, t, u)$ is substituted by 
$CI_2 (\Gamma)$ if this means $CI_2 (\Gamma; s, t, u)$ clearly.  
Similarly, for other invariants in
 this paper we admit the abbreviation like this if its meaning is  clear.  

\begin{thm}{\rm (Polyak)}
$CI_2$ is an invariant of plane curves which is as strong as the triple of the three Arnold's basic invariants 
$(J^+, J^-, St)$.  (The definitions of $J^+, J^-, St$ are
 in 
\cite{arnold1}, \cite{arnold2}, \cite{polyak1}.  )
\end{thm}

\begin{proof}
By using Polyak's formulation of the Arnold basic 
invariants \cite{polyak1}, 
the triple of the three Arnold's basic invariants $\left( J^+ (\Gamma), J^- (\Gamma),
St (\Gamma) \right)$ is equal to 
\[\left ( CI_2 (\Gamma; -\frac{1}{2}, 1, -3), 
CI_2 (\Gamma; -\frac{3}{2}, 1, -3), 
CI_2 (\Gamma; \frac{1}{4}, -\frac{1}{2}, \frac{1}{2}) \right) \]   
and three vectors $(-\frac{1}{2}, 1, -3), (-\frac{3}{2}, 1, -3), (\frac{1}{4}, -\frac{1}{2}, \frac{1}{2})$ 
are linearly independent.  These two facts imply  this theorem.  
\end{proof}

\begin{rem}
Index is independent of the basic invariants 
(cf. Figure \ref{index-cli2}).  
\begin{figure}[htbp]
\begin{center}
\includegraphics[width=15cm,height=2cm,keepaspectratio]
{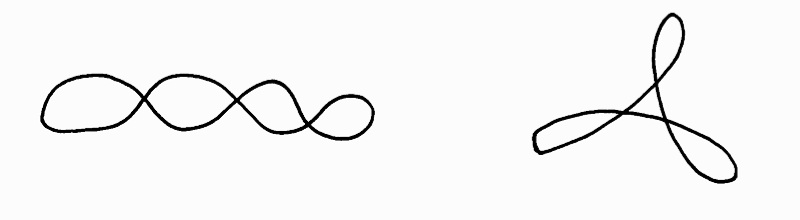}
\caption{$CI_2$ is $3s+\frac{3}{2}t$ on
 both curves, but the value of index 
is $0$ for the left figure, $2$ for the
 right figure.  }
\label{index-cli2}
\end{center}
\end{figure}
\end{rem}

We represent \[XYXYZZ-YXYZZX+XYZZXY-YZZXYX+ZZXYXY-ZXYXYZ, \]
as $[XYXYZZ]$
and
represent $XXYYZZ-XYYZZX$ as $[XXYYZZ].  $   

For every curve $\Gamma$, let $i$ be index, 
and we define $CI_3$ by 

\[CI_3 (\Gamma;  s, t) := \left[ s [XYXYZZ]+t [XXYYZZ], w_\Gamma \right]
+i.  \]

\begin{rem}
$CI_3 (\Gamma;  s, t)
 = \overline{CI}_3 (\Gamma;  s, t) 
 + i
.  $
\end{rem}

\begin{thm}\label{cli3.thm}
$CI_3$ is an invariant of plane closed curves.  
\end{thm}

\begin{proof} 
To prove this, we must prove that $CI_3 (\Gamma)$ is not varied by 
an arbitrary base point move
 (Figure \ref{basepoint-move}) for every closed curve $\Gamma$, but it is immediately concluded
 by Proposition \ref{[]}.  
\end{proof}

\begin{figure}[htbp]
\begin{center}
\includegraphics[width=20cm,height=2.5cm,keepaspectratio]
{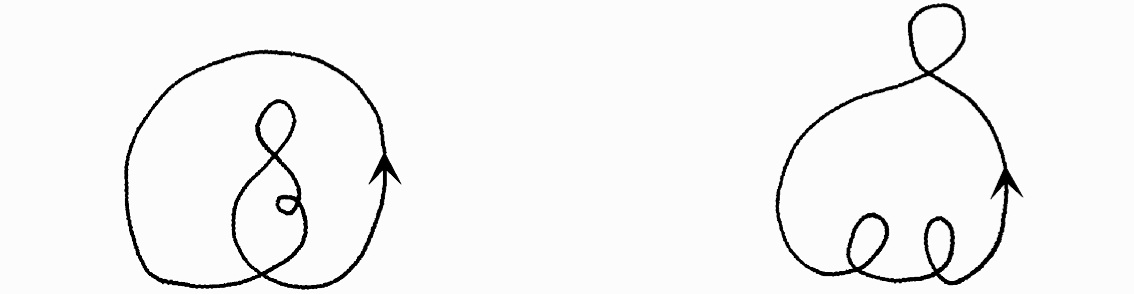}
\caption{Two curves such that the 
value of the {HOMFLY}
 polynomial is $x^2z_2$ and $(i, J^+, J^-, St)$ 
is $(2, -4, -7, 2)$. The value of $CI_3$ is 
$t+2$ for the left figure, $-t+2$ for the right figure.  }
\label{closed-curve4}
\end{center}
\end{figure}

There exist two curves such that  
the values of the HOMFLY polynomial
 of immersed plane curves \cite{CGM}, index, basic invariants are the same;  however, the value of $CI_3$ on one curve is different from that on the other  (Figure \ref{closed-curve4}).  

This example implies the following.  

\begin{cor}
$CI_3$ is independent of index, the Arnold's basic invariants and the HOMFLY polynomial
 of immersed plane curves.  
\end{cor}

As can be seen from the examples
 above (the case of $CI_2$ and $CI_3$), we can get some invariants by the  normalization
 of $\overline{CI}_n$.  
We denote by $CI_n$ a normalized 
invariants which is made of $\overline{CI}_n$.  

\begin{cor}
Suppose $-\Gamma$ has only difference of
orientation from $\Gamma$, 
and let $\Gamma^{r}$ be the reflection
 of $\Gamma$, 
\[CI_3 (-\Gamma) = -CI_3 (\Gamma), \quad CI_3
(\Gamma^r) = -CI_3 (\Gamma).  \]
\end{cor}
\begin{rem}
\[CI_2 (-\Gamma) = CI_2 (\Gamma).  \]
\end{rem}

\subsection{Strengthening $CI_n$}\label{strong-c}

We can strengthen $CI_n$ by the
 similar method of $GI_n$ in Sect. \ref{Gin}.  
We must define {\itshape marked cyclic equivalent} which is the canonical generalization of 
cyclic equivalent defined in Sect. \ref{const}.

\begin{dfn} (marked cyclic equivalent)\label{marked-jyunkaidouti}
Let $x$ be $w (2)$ $\cdots$ $w (m)$ of $w$ $=$ $w (1)$$w (2)$ $\cdots$ $w (m)$, 
for two arbitrary Gauss words $w, w' \in W_n$, relation $\sim$ is defined by
\[
w \sim w' \stackrel{\mathrm{def}}{\Longleftrightarrow} 
\left \{
\begin{array}{ll}
w=Ax~\text{and}~w'=-xA & \text{if}~d (A)=1 \\
w=\dot{A}x~\text{and}~w'=x\dot{A} & \text{if}~d (A)=2
\end{array}\right. \]  
This relation $\sim$ is called {\itshape marked cyclic equivalent}.  
\end{dfn}

Replacing cyclic equivalent by {\itshape marked cyclic equivalent}, we can easily check that the similar results are established and can easily generalize Sect. \ref{const}.     
Therefore we only see the case of $CI_3$.

For every curve $\Gamma$, let $i$ be index, 
and we define $GCI_3$ by   

$GCI_3 (\Gamma;  s, t, u) := \left[ s[XYXYZZ]+t[XXYYZZ]
+u[\dot{X}\dot{X}YY], w_\Gamma \right]
+i.  $
 
\begin{rem}
$GCI_3 (\Gamma;  s, t, u)
 = CI_3 (\Gamma;  s, t) 
 + \left[ u[\dot{X}\dot{X}YY], w_\Gamma \right]
.  $
\end{rem}
 
\begin{cor}
$GCI_3$ is an invariant of plane curves.  
\end{cor}

\begin{ex} 
There exist two curves 
such that  the value of $CI_3$, index, basic invariants are the same; however, 
the value of $GCI_3$ on one curve is
 different from that on the other (Figure \ref{closed-curve3}).

\begin{figure}[htbp]
\begin{center}
\includegraphics[width=20cm,height=3cm,keepaspectratio]
{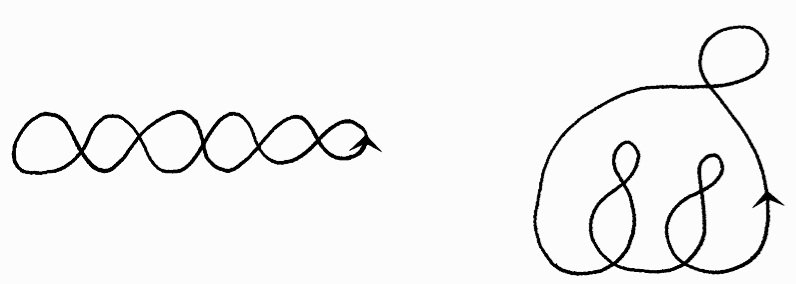}
\caption{Two curves such that 
$($$CI_3,$ $i,$ $J^+,$ $J^-,$ $St$$)$
 is $($$0,$ $0,$ $0,$ $-5,$ $0$$)$. The value of $GCI_3$ is 
$0$ for the left figure, $-8u$ for the right figure.  }
\label{closed-curve3}
\end{center}
\end{figure}
\end{ex}

In particular, this example implies the following.  

\begin{cor}
$GCI_3$ is
 a stronger invariant than $CI_3.  $
\end{cor}

\section{Application of $I_n$ to long curves}\label{long}

\subsection{Construction of invariants $LI_n$}\label{constli}

When we treated long curves via word theory, 
the following theorem is basic and fundamental.  

\begin{thm}\label{long-curve-thm}
Every long curve is represented as a nanoword.   
\end{thm}

\begin{proof}
Regard $-\infty$ on x-axis as a base point
and 
repeat the proof of Theorem \ref{curve}.  
\end{proof}

Let $i$ be index.  
By Theorem \ref{long-curve-thm}, we get 
the sequence of invariants of long curves
 $\{LI_n\}$ defined by $LI_n = I_n + i$.   

\begin{rem}
For an arbitrary function of index $i$,  $I_n + f(i)$ is plane isotopy invariant.  
\end{rem}
   
\subsection{The basic invariants of long curves}

In similar way of defining the Arnold's basic invariants 
of plane closed curve in Sect. \ref{arnold-invariant},
 we define the basic invariants of long curves 
in this subsection (cf. \cite{GN}, \cite{ZZP}).

\begin{dfn} (basic invariants of long curves)
$J^+$ is increased by 2 under positive elementary move $I\!I^+$ but not change under the other, 
$J^-$ is decrease by 2 under positive elementary move $I\!I^-$ but not change under the other, 
$St$ is increased by 1 under positive elementary move $I\!I\!I$ but not change under the other, 
and 
satisfy the following conditions 
\[J^+ (L_i)=-\left|i\right|, \quad J^- (L_i)=-2\left|i\right|, \quad
 St (L_i)=\frac{1}{2}\left|i\right| \]
for the base curves $\{L_i\}$ defined by Figure \ref{long-base}.  
\begin{figure}[htbp]
\begin{center}
\includegraphics[width=15cm,height=2cm,keepaspectratio]
{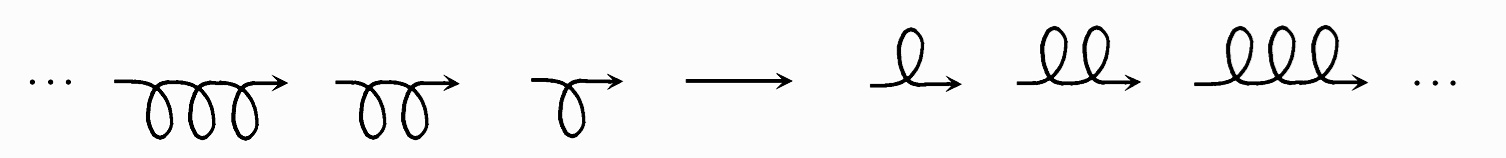}
\begin{flushleft}
\qquad \qquad \  $L_{-3}$
\qquad  \   $L_{-2}$
\qquad \ $L_{-1}$
\qquad \  $L_0$
\qquad \quad $L_1$
\qquad \   $L_2$
\qquad \quad \ $L_3$
\end{flushleft}
\caption{Base curve $\{L_i\}.  $}
\label{long-base}
\end{center}
\end{figure}
\end{dfn}

\subsection{$LI_2$ and $LI_3$}
 
For long curve L, let $i$ be index and $n$ the number of double points, we define $LI_2$ by 
\[LI_2 (L; s, t, u, v):= s n + 
\left \langle t XXYY
+u XYYX+v XYXY, w_L\right\rangle-\frac{t}{2}i^2 \]  

As case of  $CI_2$, $LI_2$ is not depend on the orientation of long curve L.

\begin{thm}\label{longthm}
$LI_2$ is an invariant of plane curves which is as strong as $(J^+, J^-, St)$.  

In other words,
 for two arbitrary long curves $L_1$, $L_2$, 
\[LI_2 (L_1)=LI_2 (L_2) \Longleftrightarrow 
 (J^+ (L_1), J^- (L_1), St (L_1))=
 (J^+ (L_2), J^- (L_2), St (L_2)).  \]
\end{thm}

Before we begin proving Theorem\ref{longthm}, 
we will prove Lemma \ref{koutousiki}.  
 (Similar formula is concluded in case closed curves \cite{polyak1}.  )
\begin{lem}\label{koutousiki}
\[LI_2(L; \frac{1}{2}, 1, 1, 1) = \frac{n}{2}+\left \langle XXYY + XYYX + XYXY, w_L \right \rangle
-\frac{i^2}{2} \equiv 0\]
In particular, left side is independent of  elementary move I\!I and I\!I\!I.    
\end{lem}
\begin{proof}
For $w_L= (\mathcal{A}, w_L)$, 
we have $i=\sum_{A \in \mathcal{A}} {\mathit{sign}}A$.    
Therefore 
\[i^2=\left  (\sum_{A \in \mathcal{A}} \mathit{sign}A \right)^2 =
\left \langle \dot{X}\dot{X} + 2XXYY + 2XYYX + 2XYXY, w_L \right \rangle.  \]
The number of double points $n$ is equal to $\left \langle \dot{X}\dot{X}, w_L \right \rangle$.  
\end{proof}

Next, we will prove
 Theorem \ref{longthm}.  
It is sufficient that we prove the following.  

\begin{proof}  ( $\Longrightarrow$).  
The following three relations are concluded by Proposition \ref{long-basic}.   
\begin{prop}\label{long-basic}
\begin{eqnarray*}\label{li2-vs-basic}
J^+ (L)&=&LI_2 (L; -\frac{1}{2}, 1, -1, -3), \\ 
J^- (L)&=&LI_2 (L; -\frac{3}{2}, 1, -1, -3), \\
St (L)&=&LI_2 (L; \frac{1}{4}, -\frac{1}{2}, \frac{1}{2}, \frac{1}{2}).  
\end{eqnarray*}
\end{prop}
 (Proof of Proposition \ref{long-basic}.  )
Let $n$ be the number of double points and
 $f(i)$ a function on index $i$.  
By using Theorem \ref{I_d}, 

 $I_2 (L; s, t, u, v)= n s +\left \langle
  t XXYY + u XYYX + v XYXY \right \rangle -t 
\frac{i^2}{2}$

is an invariant of long curves.  

By definition of $J^+, J^-, St$, 
\begin{eqnarray*}\label{x-basic}
J^+ (L)&=&LI_2 (L; s, 2s+2, 2s, 2s-2), \\ 
J^- (L)&=&LI_2 (L; s, 2s+4, 2s+2, 2s), \\
St (L)&=&LI_2 (L; s, 2s-1, 2s, 2s).  
\end{eqnarray*}

These $J^+(L), J^-(L), St(L)$ satisfies 
\[J^+ (L_i)=-\left|i\right|,  \quad J^- (L_i)=-2\left|i\right|, \quad
 St (L_i)=\frac{1}{2}\left|i\right|.  \]  

Especially, In the case $s=-\frac{1}{2}$ on $J^+$, in the case $s=-\frac{3}{2}$ on $J^-$, 
in the case $s=-\frac{1}{4}$ on St, 
these are still the basic invariants.  
We have thus proved the Proposition
 \ref{long-basic} that
 implies ($\Longleftarrow$).  

($\Longleftarrow$).  
For two arbitrary
 long curves $L_1$, $L_2$, 
assume 
\[
(J^+ (L_1), J^- (L_1), St (L_1))=
 (J^+ (L_2), J^- (L_2), St (L_2)).  \]
Four vectors $(\frac{1}{2}, 1, 1, 1)$, $(-\frac{1}{2}, 1, -1, -3)$, $(-\frac{3}{2}, 1, -1, -3)$, $(\frac{1}{4},
 -\frac{1}{2}, \frac{1}{2}, \frac{1}{2})$ 
are linearly independent.  
By adding Lemma \ref{koutousiki}
 and Proposition \ref{long-basic}, 
the assumption implies 
$LI_2(L_1)$ $=$ $LI_2(L_2)$.  
\end{proof}

\begin{rem}
Index is independent of the basic invariants.  
 (cf. Figure \ref{index-li2}).   
\begin{figure}[htbp]
\begin{center}
\includegraphics[width=10cm,height=2cm,keepaspectratio]
{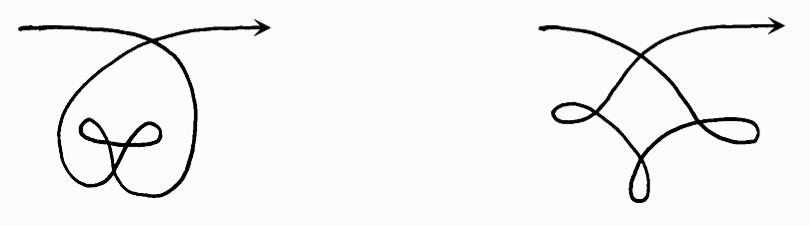}
\caption{$LI_2$ is $4s+t-3u$ on both curves,  
but 
the value of index is $0$ for the left figure,  $2$ for the right figure.  }
\label{index-li2}
\end{center}
\end{figure}
\end{rem}

\begin{rem}
If $n$ of them are not same, 
the values of $LI_2$ for two curves are not the same  because  $J^+-J^-=n$.  
\end{rem}

For every long curve $L$, let $i$ be index, 
we define $LI_3$ by 

\begin{eqnarray*}
LI_3 (L; x_1, x_2, \ldots, x_{14}, x_{15}) 
&:=& \langle x_1 XYXYZZ + x_2 XYXZZY
+x_3 XYZZXY\\
& & +x_4 XYYZXZ+x_5 XXYZYZ+x_6 XYZYZX\\
& & +x_7 XYYXZZ+x_8 XXYZZY+x_9 XYZZYX\\
& & +x_{10} XYZYXZ+x_{11} XYXZYZ+x_{12} XYZXZY\\
& & +x_{13} XYZXYZ+x_{14} XXYYZZ
+x_{15} XYYZZX, \\
& & w_L \rangle
+i.  
\end{eqnarray*}

\begin{thm}
$LI_3$ is an invariant of long curves.  
\end{thm}
\begin{proof}
Theorem \ref{I_d} immediately deduces this conclusion.  
\end{proof}

\begin{cor}
Let $L^r$ be the reflection of a long curve $L$, 
$LI_3 (L^r)=-LI_3 (L)$.  
\end{cor}
\begin{rem}
$LI_2 (L^r)=LI_2 (L) $.  
\end{rem}

The following are examples of several pairs of long curves such that the values of index and three 
basic invariants are the same, but the value of  $LI_3$ on one long curve 
is different from that on the other.   

\begin{ex}
For $L_1$, $L_1^r$, $L_2$, $L_2^r$ (Figure \ref{long-curve6})
such that  $(i, J^+, J^-, St)$ is $(0, 0, -4, 0)$, each value of  $LI_3 $ is different from another of them.

\begin{figure}[htbp]
\begin{center}
\includegraphics[width=12cm,height=2cm,keepaspectratio]
{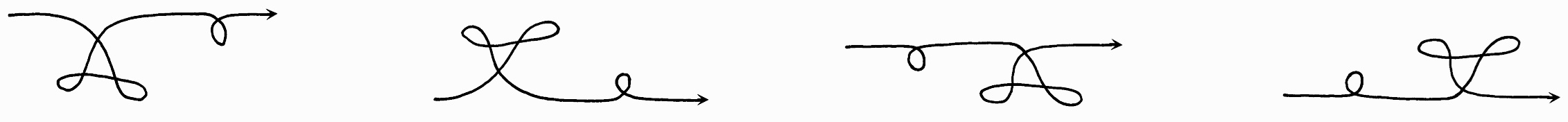}

\begin{flushleft}
\qquad\qquad\qquad
$L_1$\qquad\qquad\qquad\quad   
$L_1^r$\qquad\qquad\qquad  
$L_2$\qquad\qquad\qquad\quad  
$L_2^r$
\end{flushleft}

\caption{Four long curves such that  $(i, J^+, J^-, St)= (0, 0, -4, 0).  $
$LI_3 (L_1)$ $=$ $2 x_7$ $-$ $x_{14}$ 
$-$ $x_{15},$ 
$LI_3 (L_1^r)$ $=$ $-2 x_7$ $+$ $x_{14}$ $+$ $x_{15},$ 
$LI_3 (L_2)$ $=$ $2 x_8$ $-$ $x_{14}$
 $-$ $x_{15},$ 
$LI_3 (L_2^r)$ $=$ $-2 x_8$ $+$ $x_{14}$
 $+$ $x_{15}.$  }
\label{long-curve6}
\end{center}
\end{figure}
\end{ex}

\begin{ex}
For  $L_3$, $L_4$, $L_5$, $L_6$ (Figure \ref{long-curve7})
 such that $($$i,$ $J^+,$ $J^-,$ $St$$)$ $=$ $($$2,$ $-4,$ $-8,$ $2)$, 
each value of $LI_3$ is different from another of them.

\begin{figure}[htbp]
\begin{center}
\includegraphics[width=12cm,height=5cm,keepaspectratio]
{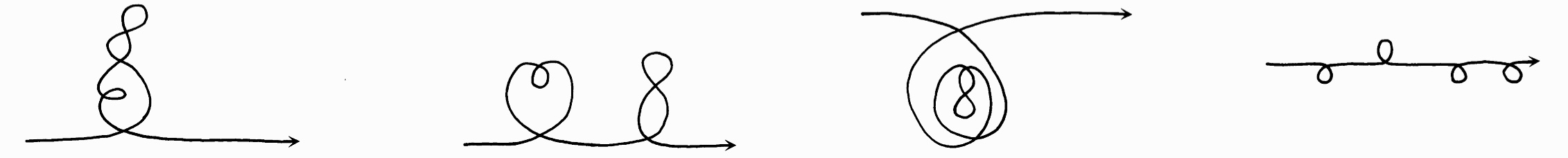}

\begin{flushleft}
\qquad\qquad\qquad
$L_3$\qquad\qquad\qquad\quad   
$L_4$\qquad\qquad\qquad  
$L_5$\qquad\qquad\qquad\quad  
$L_6$
\end{flushleft}

\caption{Four long curves such that  $ ($$i,$ $J^+,$ $J^-$ $St$$)$= $($$2,$ $-4,$ $-8,$ $2$$)$$.$  $LI_3 (L_3)$ $=$ $-x_7$ $-$ $x_9$ $+$ $2,$ 
$LI_3 (L_4)$ $=$ $-2 x_8$ $+$ $2,$ 
$LI_3 (L_5)$ $=$ $2 x_9$ $+$ $2,$ 
$LI_3 (L_6)$ $=$ $2 x_{14}$ $+$ $2.$  }
\label{long-curve7}
\end{center}
\end{figure}
\end{ex}

\subsection{Strengthening $LI_n$}\label{strong-l}

We can strengthen $LI_n$ in the
 way constructing $GI_n$ in Sect. \ref{Gin}.  

For every long curve $L$, let $i$ be index, 
we define $GLI_3$ by 

\begin{eqnarray*}
GLI_3 (L; x_1, x_2, \ldots, x_{20}, x_{21}) 
&:=& LI_3 (L; x_1, x_2, \ldots, x_{14}, x_{15}) \\
& & 
+ \langle x_{16} \dot{X}\dot{X}YY + 
x_{17} \dot{X}YY\dot{X}+x_{18} XX\dot{Y}\dot{Y}\\
& & +x_{19} X\dot{Y}\dot{Y}X+x_{20} X\dot{Y}X\dot{Y}+x_{21} \dot{X}Y\dot{X}Y, w_L \rangle.  
\end{eqnarray*}

\begin{cor}
$GLI_3$ is an invariant of long curves. 
\end{cor}

\begin{ex} 
There exist two curves 
such that  the value of $LI_3$, index, basic invariants are the same; however, 
the value of $GLI_3$ on one curve is
 different from that on the other (Figure \ref{long-curve8}).

\begin{figure}[htbp]
\begin{center}
\includegraphics[width=10cm,height=3cm,keepaspectratio]
{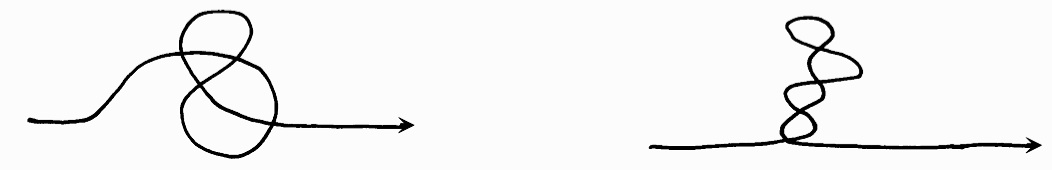}
\caption{Two curves such that 
$($$LI_3,$ $i,$ $J^+,$ $J^-,$ $St$$)$
 is $($$0,$ $0,$ $0,$ $-4,$ $0$$)$. The value of  $GLI_3$ is 
$0$ for the left figure,  $-2x_{17}$ $+$ $2x_{19}$ for the right figure.  }
\label{long-curve8}
\end{center}
\end{figure}
\end{ex}

In particular, this example implies the following.  

\begin{cor}
$GLI_3$ is a stronger invariant than $LI_3.$  
\end{cor}

\subsection{The Arnold-type invariant of degree 3}\label{degree3}

\begin{dfn}(The Arnold-type invariant of degree 3)
Let 
$J^{+} $-type invariant be an invariant do not change under elementary move $I\!I^-$, $I\!I\!I$, 
$J^{-} $-type invariant be an invariant do not change under elementary move~ $I\!I^+$, $I\!I\!I$, 
$St$-type invariant be an invariant do not change under elementary move~ $I\!I^+, I\!I^-$.  
\end{dfn}

$J^{+} $-type invariant, $St $-type invariant are available from  $GLI_3$.  
These are expressed  $J^+_3$ ($J^+$ of degree 3),  $St_3 $ ($St $ of degree 3).  
  
\begin{thm}Let L be long curve.  
\begin{eqnarray*}
J^+_3 (L)&=&GLI_3 (L; s+t-v, s-t+u, -s+t+v, -s+3t-u, s+t-v, \\
& & 2t-v, s, 
s, t, u, 2t-v, -2s+4t-u, -2s+2t+v, s-t+v, \\
& & \frac{1}{2}t, \frac{1}{2}s, 
v, \frac{1}{2}s, \frac{1}{2}t, -s+2t-\frac{1}{2}u, \frac{1}{2}u), \\
St_3 (L)&=&GLI_3 (L; 2u, p, 2s, q, 2x, 2y,
 2u, 2x, 2y, 2z, r, 
2z, 2z, s, t,
 u, v, x, y, z, z).  \\  
\end{eqnarray*}
\end{thm}

\begin{proof}
To prove this theorem, we check 
variations of each parameter's coefficient 
of $GLI_3$: 
\begin{eqnarray*}
GLI_3 (L; x_1, x_2, \ldots, x_{20}, x_{21}) 
&:=& LI_3 (L; x_1, x_2, \ldots, x_{14}, x_{15})\\
& & 
+ \langle x_{16} \dot{X}\dot{X}YY + 
x_{17} \dot{X}YY\dot{X}+x_{18} XX\dot{Y}\dot{Y}\\
& & +x_{19} X\dot{Y}\dot{Y}X+x_{20} X\dot{Y}X\dot{Y}+x_{21} \dot{X}Y\dot{X}Y, w_L \rangle.  
\end{eqnarray*}

Consider variations of each parameter's coefficient under each elementary move.  
\end{proof}

\section{Application of $I_n$ to fronts}\label{front}
\subsection{The basics of fronts}\label{dfn-front}

To begin with, we recall basic concepts or results about fronts.    

\begin{figure}[htbp]
\begin{center}
\includegraphics[width=10cm,height=3cm,keepaspectratio]
{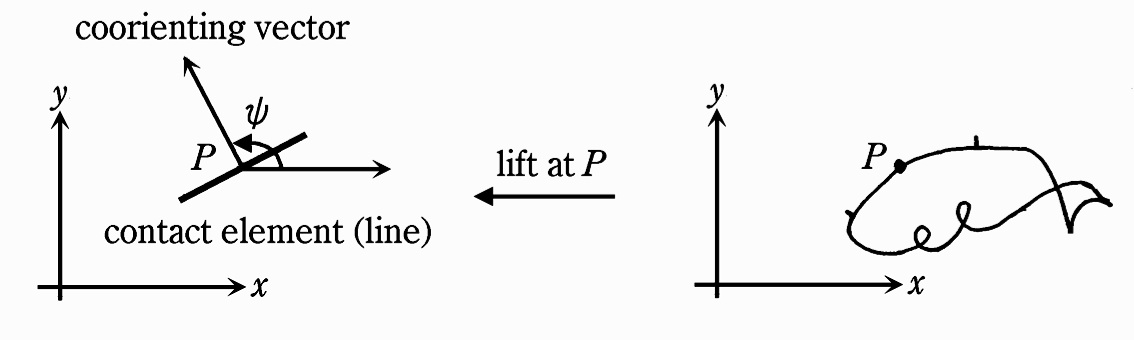}
\caption{Contact element.  }
\label{contact-element}
\end{center}
\end{figure}

\begin{dfn}(contact element)
A {\itshape contact element} in the plane is a line in the tangent plane (Figure \ref{contact-element}).  
\end{dfn}

\begin{dfn}(coorientation)
The {\itshape coorientation}  of a contact element 
is the choice of the half-plane into which the contact element divides the tangent plane  (Figure \ref{contact-element}).  
\end{dfn}

The manifold $M$ consisting of all of the contact elements in the plane is diffeomorphic to solid torus.  
Consider immersion : $S^1 \to M$.  
For each element of this curve in the plane, the coorient of the element is determined  (Figure \ref{contact-element}).  
This curve which is Legendrian submanifold of  $M$ 
is called Legendrian curve, 
the image of projection to plane of this curve is 
called a {\itshape front}.  

\begin{figure}[htbp]
\begin{center}
\quad \ $2k$ cusps \qquad\qquad\qquad \quad
$i-1$ kinks \ \  $2k$ cusps
\includegraphics[width=10cm,height=3cm,keepaspectratio]
{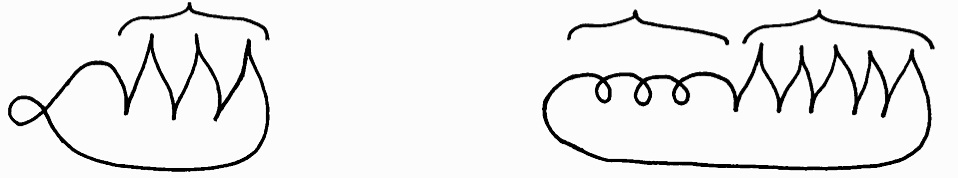}

$K_{0, k}$ \qquad\qquad\qquad
\qquad\qquad\qquad \quad $K_{i, k}$
\caption{Base fronts $\{K_{i, k}\}.  $}
\label{front-base}
\end{center}
\end{figure}

By using the following several concepts: elementary move  (Definition \ref{f-kihon}), index  (Definition \ref{index}), and Maslov index  (Definition \ref{maslov}), 
a front is regarded as a plane curve generated by $K_{i, k}$  (Figure \ref{front-base}) via elementary 
moves because Theorem \ref{gr-thm} is established.

Let the alphabet $\alpha_* $, the  involutions $\tau_1 $, and $\tau_2 $,
 the set $S_{*} $ be 
\begin{eqnarray*}
\alpha_{*} &=& \{a_+, a_-, b_+, b_-\}, \\
\tau_1& : &a_{+}\mapsto b_{+}, a_{-}\mapsto b_{-}, \\
\tau_2& : &a_{+}\mapsto b_{-}, a_{-}\mapsto b_{+}, \\
S_{*} &=& \{ (a_{\epsilon}, a_{\eta}, a_{\xi}),
  (b_{\epsilon}, b_{\eta}, b_{\xi}) :\epsilon, \eta, \xi = +,~{\rm or}~-\}.  
\end{eqnarray*}

 (They are different from the involution  $\tau $,  $S_{*} $ for  $\alpha_{*} $ in \cite{turaev2}.  )

To the following relations
\[a_+ + b_+ =0, \quad a_- + b_- =0\]
are established, 
we consider the ring 
 $\overline{\mathbf{Q}\tilde{\alpha_{*}}} $ 
into which two-sided ideal generated by
  $a_+ + b_+ $ ,   $a_- + b_-$ divides 
 the monoid algebra $\mathbf{Q}\tilde{\alpha_{*}}$
 where $\tilde{\alpha_{*}}$ is commutative monoid generated
 by $\alpha_{*} \cup \{1\}$.

\begin{rem}
If we regard  $a_{+} $,  $b_{+} $,  $\tau_1 $ in the same light as  $a $,  $b $,  $\tau $ the discussion from Sect. \ref{id} to Sect. \ref{long} still holds
on this establishment 
because  $a_{-} $,  $b_{-} $, $\tau_2 $ do not 
appear Sect. \ref{id}, \ref{plane}, and  \ref{long}.  
That is to say, this establishment is canonical generalization of the establishment in Sect. \ref{id}, \ref{plane}, and \ref{long}.  
\end{rem}

From this section, we denote by $|~|$
 projection $:$
 $\mathcal{A} \to \alpha_{*}.  $

\begin{dfn} ($|A|$ corresponding to a double point)
We define $|A|$ corresponding to a double point $A$ by 
Figure \ref{co-orient2}.  
\begin{figure}[htbp]
\begin{center}
\includegraphics[width=10cm,height=3cm,keepaspectratio]
{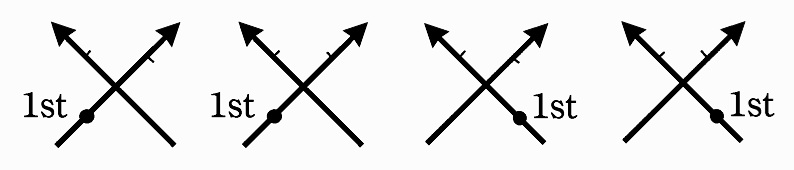}

{\Large $a_+$ 
\qquad\ \ $a_-$
\qquad\ \ $b_+$
\qquad\ \ \ $b_-$}
\caption{Correspondence 4 patterns of double points to 4 alphabets.  }
\label{co-orient2}
\end{center}
\end{figure}
\end{dfn}

\begin{dfn} ($|K|$ corresponding to a cusp)
$|K|$ where  $K$ is a cusp is defined by 
Figure \ref{cusps1}.  
\begin{figure}[htbp]
\begin{center}
\includegraphics[width=10cm,height=3cm,keepaspectratio]
{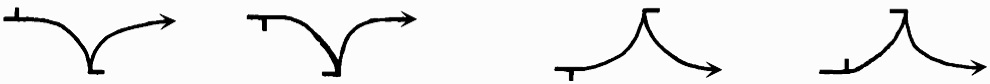}

{\Large \ $a_+$ 
\qquad\ \ \  $a_-$
\qquad\qquad \ $b_+$
\qquad\ \ \ \ $b_-$}
\caption{Correspondence 4 patterns of cusps to 4 alphabets.  }
\label{cusps1}
\end{center}
\end{figure}
\end{dfn}

\begin{dfn}(index)\label{index}
{\itshape Index} of a front $F$  is the number of the full rotations of the coorienting normal vector counter clockwise while we trip along the front once.  
\end{dfn}

\begin{dfn}(Maslov index)\label{maslov} 
For \'{e}tale word  $ (\mathcal{A}_F,w_F) $ of front $F$, let $K$ is a letter of $\mathcal{A}_F$
 for a cusp.  
Maslov index  $\mu $ is defined by \[\mu=card\{K \in \mathcal{A}_F | |A|=a_+ ~{\rm or}~ b_+\}-
card\{K \in \mathcal
{A}_F | |A|=a_- ~{\rm or}~ b_-\}.  \]    
\end{dfn}

The local moves of fronts (Figure \ref{safe-move}, Figure \ref{dangerous-move}, Figure \ref{direct-cusp-crossing1}, 
Figure \ref{inverse-cusp-crossing}, Figure \ref{cusp-birth1}) is defined by the  following.  
Suppose the local moves are admitted apart from base point.  

\begin{dfn}\label{f-kihon} (elementary move of fronts)
Let $x$, $y$, $z$ be words that consist of the letter in $\mathcal{A}$.  Four kinds of elementary move  $I\!I$ are defined by
\begin{eqnarray*} 
S\!I\!I^+ : (\mathcal{A}\cup\{A, B\}, xAByABz) \to  (\mathcal{A}, xyz) ~\text{with}~ \tau_1 (|A|)=|B|, \\
S\!I\!I^- : (\mathcal{A}, xyz) \to  (\mathcal{A}\cup\{A, B\}, xAByBAz) ~\text{with}~ \tau_1 (|A|)=|B|, \\
D\!I\!I^+ :  (\mathcal{A}, xyz) \to  (\mathcal{A}\cup\{A, B\}, xAByABz) ~\text{with}~ \tau_1 (|A|)=|B|,  \\
D\!I\!I^- :  (\mathcal{A}\cup\{A, B\}, xAByBAz) \to  (\mathcal{A}, xyz) ~\text{with}~ \tau_1 (|A|)=|B|.  
\end{eqnarray*}

Let $x$, $y$, $z$, t be words that consist of the letter in $\mathcal{A}-\{A,B,C\}$.  

elementary move $I\!I\!I$ is defined by 

\[I\!I\!I :  (\mathcal{A} , xAByACzBCt) \to 
 (\mathcal{A}, xBAyCAzCBt)
 ~\text{for}~ (|A|, |B|, |C|)\in 
S_{*}.  \] 

Let $x$, $y$, $z$ be words that consist of the letter in $\mathcal{A}-\{K\}$.  

Two kinds of elementary move $\Pi$ is defined by
\begin{eqnarray*}
\qquad\quad \Pi^+ :  (\mathcal{A}, xKyz) \to
  (\mathcal{A}\cup\{A, B\}, xAKByABz)
\qquad \qquad \qquad \qquad \qquad \qquad \qquad\\ 
\text{or}~ 
(\mathcal{A}, xKyz) \to  (\mathcal{A}\cup\{A, B\}, xAByAKBz)~\text{with}~ \tau_2 (|A|)=|B|, 
\qquad\\
\qquad\quad \Pi^- :  (\mathcal{A}\cup\{A, B\}, xAKByBAz)  \to  (\mathcal{A}, xKyz)
\qquad \qquad \qquad \qquad \qquad \qquad \qquad\\ 
\text{or}~ 
 (\mathcal{A}\cup\{A, B\}, xAByBKAz) \to  (\mathcal{A}, xKyz)~\text{with}~ \tau_2 (|A|)=|B|.  
\end{eqnarray*}
 
Let $x$, $y$ be words that consist of the letter in $\mathcal{A}.$  elementary move $\Lambda$ is 
defined by 

\[\Lambda :  (\mathcal{A}, xy) \to  (\mathcal{A}\cup\{A, K_1, K_2\}, xAK_1K_2Ay) 
~\text{with}~ \tau_1\circ\tau_2 (|K_1|)=|K_2|.  \]  

The positive elementary moves is the above direction, negative is the
 inverse direction.

\begin{figure}[htbp]
\begin{center}
\includegraphics[width=12cm,height=3cm,keepaspectratio]
{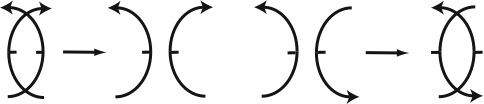}
\caption{Elementary move   $S\!I\!I^+ $ (left figure),  $S\!I\!I^-$ (right figure).  }
\label{safe-move}
\end{center}
\end{figure}

\begin{figure}[htbp]
\begin{center}
\includegraphics[width=12cm,height=3cm,keepaspectratio]
{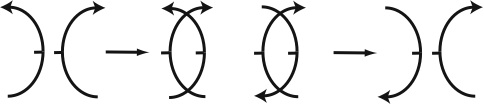}
\caption{Elementary move   $D\!I\!I^+ $ (left figure),  $D\!I\!I^- $ (right figure).  }
\label{dangerous-move}
\end{center}
\end{figure}

\begin{figure}[htbp]
\begin{center}
\includegraphics[width=12cm,height=3cm,keepaspectratio]
{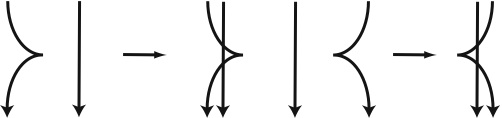}
\caption{Elementary move   $\Pi^+ $.  }
\label{direct-cusp-crossing1}
\end{center}
\end{figure}
\begin{figure}[htbp]
\begin{center}
\includegraphics[width=12cm,height=3cm,keepaspectratio]
{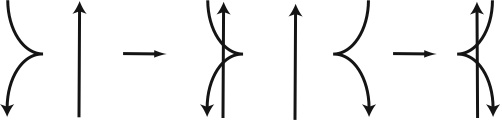}
\caption{Elementary move   $\Pi^-$.  }
\label{inverse-cusp-crossing}
\end{center}
\end{figure}
\begin{figure}[htbp]
\begin{center}
\includegraphics[width=12cm,height=3cm,keepaspectratio]
{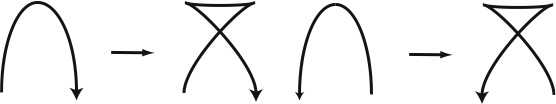}
\caption{Elementary move $\Lambda$.  }
\label{cusp-birth1}
\end{center}
\end{figure}
\end{dfn}

\begin{rem}
$S\!I\!I$ stands for safe 2-move,and $D\!I\!I$ means dangerous 2-move  (2-move, 3-move resemble Reidemeister move $I\!I$, $I\!I\!I$).  
In the lift of the plane, self-tangency under dangrous move is corresponded  to 
crossing of the Legendrian curve
 (cf. \cite{aicardi2}).  
\end{rem}

We have the next theorem 
due to Gromov \cite{gromov}.  

\begin{thm}{\rm(Gromov)}\label{gr-thm}
Any front whose index is $i$, Maslov index is $2 k$ is 
deformed via $D\!I\!I^+$, $D\!I\!I^-$,
 $S\!I\!I^+$, $S\!I\!I^-$, $I\!I\!I$, $\Pi^+$, $\Pi^-$, $\Lambda$ from $K_{i, k} $ (Figure \ref{front-base}).    
\end{thm}

The object of this section is giving a 
classification of fronts 
more detail than the classification by Theorem \ref{gr-thm}.  
To do this, we consider an application a method like 
Sect. \ref{plane}, \ref{long} to front which may 
have not only double points but also cusps.  
We consider only double points and cusps due to Theorem \ref{gr-thm}.

We get following theorem.    

\begin{thm}
All fronts are represented as 
\'{e}tale words.    
\end{thm}

\begin{proof} For an arbitrary given front  which has precisely  $n$ double points and precisely $m$ cusps, 
we name double points $A_1, A_2, \ldots, A_n$ 
and name cusps $K_1, K_2, \ldots, K_m$
along front from the base point.  
Every double point and every cusp precisely corresponds to a unique element  of $\alpha_{*}$ in Figure \ref{co-orient2}.  
\end{proof}

\subsection{The basic invariants of fronts}

The basic invariants of fronts  $J^+, J^-, St $ are 
three invariants can be defined as follows 
 (cf. \cite{arnold2}, \cite{polyak1}).  

\begin{dfn}
$J^+$ is increased by $2$ under positive elementary move  $D\!I\!I^+$, $D\!I\!I^-$ but not change under the other, 

$J^-$ is decreased by $2$ under positive elementary move  $S\!I\!I^+$, $S\!I\!I^-$ but not change under the other,

 $St$ is increased by $1$ under positive elementary move $I\!I\!I$, increased by $\frac{1}{2}$ under positive elementary move $\Pi^+$, and 
decreased by $\frac{1}{2}$ under positive elementary move $\Pi^-$,  
but not change under the other, 
and 
satisfy the following conditions
\[J^+ (K_{0, k})=-k, \quad J^- (K_{0, k})=-1, \quad St (K_{0, k})=\frac{k}{2}, \]
\[J^+ (K_{i+1, k})=-2i-k, \quad J^- (K_{i+1, k})=-3i, \quad St (K_{i+1, k})=i+\frac{k}{2} \]
for the base curves  $\{K_i\} $ defined by Figure \ref{front-base}.  

\end{dfn}

\subsection{Construction of invariants $FI_n$}\label{constfi}

In this subsection 
we will construct invariants of fronts.

\begin{dfn}(fake Gauss word)\label{fake-Gauss-word}
We call a word $w$ a {\itshape fake Gauss word 
of dimension $n$}
if $w' : 2m \to \mathcal{A}$ is a Gauss word where $w'$ is a sub-word 
of $w$ and $\{w(1), w(2),  \ldots, w(2m)\}$ $\sqcup$ $\{K_1, K_2, \ldots, K_{2l}\}$ $=$ $\{w(1), w(2),  \ldots, w(2(m+l))\}$ where $n=m+2l$.  
\end{dfn}

\begin{dfn}(fake nanoword)\label{fake-nanoword}
An \'{e}tale word $(\mathcal{A}, w)$
 over $\alpha_{*}$ is a {\itshape fake nanoword over 
$\alpha_{*}$ of dimension $n$
} if $(\mathcal{A}, w)$ satisfies $(\mathcal{A}-\{K_1, K_2, \ldots, K_{2l}\}, w')$ is an nanoword over $\alpha_{*}$ where 
the length of $w$ is 
$card\mathcal{A}' + 2l = n$ and
 $w'$ is a sub-word of $w$.  
\end{dfn}

\begin{rem}
We can consider a front $F$ on a surface which is 
a curve on a surface with coorientation and 
cusps.  
A fake nanoword gives rise to a nanoword
 $w$ by neglecting letters for cusps.  
We can calculate the genus $g(F)$ of a surface on 
which a fronts because $g(F)$ 
is equal to $g(\Gamma)$ where $\Gamma$ is
 determined by the nanoword $w$ by using Theorem
 \ref{word-curve}.  
\end{rem}

\begin{dfn}(fake word space)\label{fake-word-space}
The fake word space of degree $n$ is the word space generated by
 fake Gauss words of dimension $n$.  
\end{dfn}

In \ref{constfi} and \ref{f2-f3}, suppose fake Gauss word, fake nanoword and fake word space are made of only $1$-dimensional letters.  We denote by $FW_n$ the fake word space is made of fake Gauss words such that all letters are $1$-dimensional letters.  

Let $\epsilon$ be $\epsilon = +, ~\text{or}~-$.  
For every $A \in \alpha_{*}$-alphabet $\mathcal{A}$, 
 $\mathit{sign} (|A|)$ is  $1$ if $|A|=b_{\epsilon}$ and
 $\mathit{sign} (|A|)$ is $-1$
if $|A|=a_{\epsilon}$ (We consider its generalization
 in Sect. \ref{strong}).  
Let  $w_F$ stand for a word which is determined
by a front $F$.  
By using Theorem \ref{I_d}, 
for every pointed front (which means every front with a base point), 
the following $\{FI_n\}$ is plane isotopy invariants sequence.  
\[FI_n (F)=\left\langle \sum_k x_k v_k, w_F \right\rangle
\qquad  (n \in \mathbf{N}).  \]
${v_k}$ is a sequence consisting of all elements in  $FW_n$, and 
each $x_k$ is a parameter.

In the same way of Sect. \ref{const}, in order to simplify following description, 
we define proper equivalent classes of
 fake nanowords over $\alpha_{*}$.

\begin{dfn} (cyclic equivalent for fronts)\label{f-jyunkaidouti}
Let  $w (2) \cdots w (n)$ of $w=w (1)w (2) \cdots w (n)$ 
represent  $x$, 
for two arbitrary $w, w' \in FW_n$, relation $\sim$ is defined by
\[
w \sim w' \stackrel{def}{\Longleftrightarrow} 
\left \{
\begin{array}{ll}
w=Ax ~\text{and}~ w'=-xA & 
\text{if~A means a double point}\\
w=Kx ~\text{and}~ w'=xK & \text{if~when K means a cusp}\\
\end{array}\right. \] 
This relation  $\sim$ is called the
 {\itshape cyclic equivalent for fronts}.  
\end{dfn}

\begin{cor}
The cyclic equivalent for fronts is equivalent relation.  
\end{cor}

Let $\overline{FW}_n$ be a module consisting of  cyclic equivalent classes for front and 
$FN_d$ the $\mathbf{Q}$-module generated by 
the set of fake nanoword over $\alpha_{*}$ 
$\{(\mathcal{A}, w)\}$ such that $card\mathcal{A}$
 is less than $d+1$.  

The mapping $\left[, \right] : \overline{FW}_n \times FN_d \to \mathbf{Q}$
is defined by
\[\left[[v], w \right]:=
\left<v_1+v_2+ \cdots +v_{2l} , w_F \right>, \] 
$v_1, v_2, \ldots, v_{2l} $consist of all elements of $[v]$.

\begin{prop}\label{[f]}
For every $[w]$ which is one of base of
$\overline{FW}_n$, 
$[[w], w_F]$ is an invariant of curve  $\Gamma$.    
\end{prop}

\begin{proof}
The proof is similar to the proof of
 Proposition \ref{[]}.  
\end{proof}

\subsection{$FI_2$ and $FI_3$}\label{f2-f3}
For every front $F= (\mathcal{A}_F, w_F) $, 
let $i$ be index and $2c$ the number of cusps, 
\[n_{\epsilon}=card\{A \in \mathcal{A}_F:|A|=a_{\epsilon}~
\text{or}~b_{\epsilon}\}~ (\epsilon = +,  \text{or}~-).  \]
We define $FI_2$ by 

\begin{eqnarray*}
FI_2 (F; p, q, r, s, t, u, v) &:=& p n_+ + qn_- + \langle
 r XXYY - r XYYX + s XYXY \\
& & + t KXX + t XXK - t XKX + u KK,
 w_F \rangle +vc \\
& & +\frac{r}{2}-\frac{r}{2} i^2.  
\end{eqnarray*}

$FI_2$ is also not depend on the orientation of each front 
$F$ because this formula is symmetric.

\begin{thm}
$FI_2$ is an invariant which is stronger than$(J^+, J^-, St)$.  

$i. e. $ For a front $F$, $FI_2 (F)$ is not depend on the choice of
 the base point, 
for two arbitrary fronts $F_1$, $F_2$, 
\[FI_2 (F_1)=FI_2 (F_2) \Longrightarrow 
 (J^+ (F_1), J^- (F_1), St (F_1))=
 (J^+ (F_2), J^- (F_2), St (F_2)), \]
and the converse can not be established.  
\end{thm}
\begin{proof}
 (I) First, we will prove that for an arbitrary front $F$, $FI_2 (F)$ is not depend on the choice of a base point.  
Base point move  (Figure \ref{basepoint-move-front}) means 
 $x, y $ is word consisting of letters in $\mathcal{A}-\{A\} $,  
in the case of $A$ is a double point, as Proposition \ref{[]}, 
 $ (\mathcal{A}, w_F)=AxAy 
\to xA_{\tau} yA_{\tau}, |A_\tau|=\tau_{1} (|A|),  $
in the case of  $A $ is a cusp  $i. e. A=K $, 
 $ (\mathcal{A}, w_F)=Kxy 
\to xyK, |K|$ is not change .

$\left\langle XXYY-XYYX, w_F\right\rangle,  
\left\langle XYXY, w_F\right\rangle$ is not depend 
on the choice of a base point due to 
\cite{polyak1}.

Therefore to prove it is checking increase and decrease of only the terms 
 \[\left\langle tKXX+tXXK-tXKX, w_F\right\rangle \]
under an arbitrary base point move.

 (I\!I) For front $F_1$, $F_2$, 
\[FI_2 (F_1)=FI_2 (F_2) \Longrightarrow 
 (J^+ (F_1), J^- (F_1), St (F_1))=
 (J^+ (F_2), J^- (F_2), St (F_2)), \]
and the converse can not be established.  
 
($\Longrightarrow $).  
The following three relations are concluded by Polyak in \cite{polyak1}.   
\begin{eqnarray*}
J^+ (F)&=&FI_2 (F; -\frac{1}{2}, -\frac{3}{2}, 1, -3, \frac{1}{2}, \frac{1}{4}, -\frac{3}{4}), \\
J^- (F)&=&FI_2 (F; -\frac{3}{2}, -\frac{1}{2}, 1, -3, \frac{1}{2}, \frac{1}{4}, \frac{1}{4}), \\
St (F)&=&FI_2 (F; \frac{1}{4}, \frac{1}{4}, -\frac{1}{2}, \frac{1}{2}, -\frac{1}{4},
 -\frac{1}{8}, \frac{3}{8}).  
\end{eqnarray*}
 ( $\Longleftarrow $ can not be established.  )
There exists counterexample : Figure \ref{front1}.  

\begin{figure}[htbp]
\begin{center}
\includegraphics[width=12cm,height=3cm,keepaspectratio]
{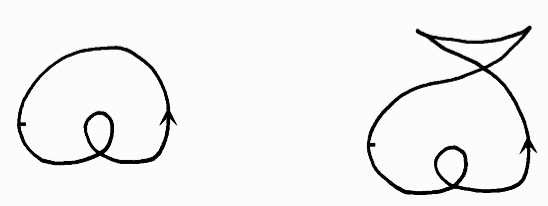}
\caption{The values of basic invariants of two fronts are 
same, but the value of  $FI_2$ of one of them  is different from that of another.  }
\label{front1}
\end{center}
\end{figure}
\end{proof}

\begin{figure}[htbp]
\begin{center}
\includegraphics[width=14cm,height=3cm,keepaspectratio]
{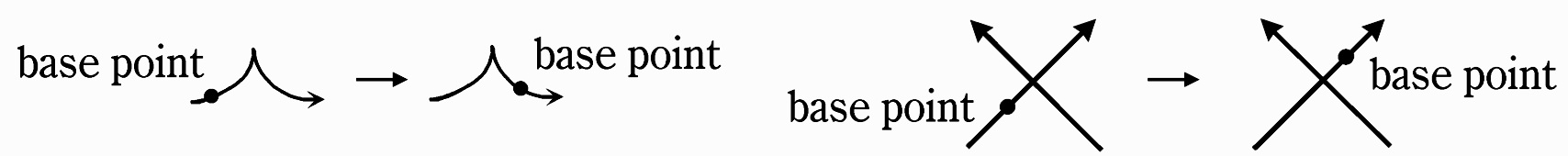}
\caption{Base point move.  }
\label{basepoint-move-front}
\end{center}
\end{figure}

\begin{rem}
There is a relation  $J^+-J^-=n_+ - n_- -c $ (cf. \cite{polyak1}).  
\end{rem}

We will consider application  $I_3 $ to fronts.  

Let $i$ be index.  
For every front $F$, 
we define $FI_3$ by 

$
FI_3 (F; x, y, z, p, q, r, s, t)
 := CI_3 (F; x, y, z)
+\Big{[}p[XKXYY]+q[KXXYY]\newline
+r[XKYXY] 
+s[XXKK] 
+t[KKK], w_F \Big{]}.  $

\begin{thm}
 $FI_3 $ is an invariant of fronts.  
\end{thm}

\begin{proof}
To prove this, we must prove that  $FI_3 (F) $ is independent of base point move
 (Figure \ref{basepoint-move-front}) for  every front, but
it is immediately concluded by
Proposition \ref{[f]}.  
\end{proof}

\begin{cor}
Suppose $-F$ has only difference of orientation from $F$, 
and let  $F^{r}$ be the reflection
 of $F$, 
\[FI_3 (-F)=-FI_3 (F), \quad FI_3 (F^r)=-FI_3 (F).  \]
\end{cor}
\begin{rem}
$FI_2 (-F)=FI_2 (F)$.  
\end{rem}

The following is example of the pair of curves such that the values of index, 
Maslov index, the basic invariants, and  $FI_2$ are the same, but the value of $FI_3$ on one front 
is different from that on the other.    

\begin{ex}
There exist two curves (Figure \ref{front2}), 
such that $($$i,$ $\mu,$ $J^+,$ $J^-,$
 $St,$ $FI_2$$)$
 $=$ $($$1,$ $2,$ $-3,$ $-4,$ $6,$
 $2 p$ $-$ $r$ $-$ $u$ $+$ $v$$)$ and the values of one 
front is different from another.  
\begin{figure}[htbp]
\begin{center}
\includegraphics[width=15cm,height=3cm,keepaspectratio]
{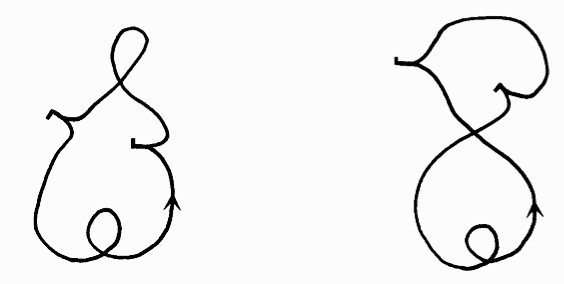}
\caption{Two fronts such that 
$($$i,$ $\mu,$ $J^+,$ $J^-,$ $St,$
 $FI_2$$)$
 $=$ $($$1,$ $2,$ $-3,$ $-4,$ $6,$
 $2 p$ $-$ $r$ $-$ $u$ $+$ $v$$)$.  
The value of $FI_3$ is 
 $2s + 1$ for left 
figure,  $2s-4p+1$ for right figure.  }
\label{front2}
\end{center}
\end{figure}
\end{ex}

\subsection{Strengthening $FI_n$}

In this section, suppose fake Gauss word  (Definition \ref{fake-Gauss-word}), fake nanoword (Definition \ref{fake-nanoword}), and fake word space (Definition \ref{fake-word-space}) are made of not only $1$-dimensional letters but also $2$-dimensional letters.  
We can strengthen $FI_n$ by the
 similar method of Sect. \ref{Gin}, 
\ref{strong-c}, \ref{strong-l}.  
We must define {\itshape marked cyclic equivalent for fronts} which is the canonical generalization of 
cyclic equivalent for fronts defined in Sect. \ref{constfi}.  

In distinction from $FW_n$, 
we denote by $\mathcal{FW}_n$
 the fake word space may have fake Gauss words which letters are not only $1$-dimensional letters but also $2$-dimensional letters.  

\begin{dfn} (marked cyclic equivalent for fronts)\label{f-jyunkaidouti}
Let  $w (2) \cdots w (n) $ of $w=w (1)w (2) \cdots w (n)$ 
represent $y$, 
for two arbitrary $w, w' \in \mathcal{FW}_n$, relation $\sim$ is defined by
\[
w \sim w' \stackrel{def}{\Longleftrightarrow} 
\left \{
\begin{array}{lll}
\text{when X means a double point,} & w=Xy ~\text{and}~ w'=-yX & 
\text{if}~d (X)=1\\
& w=\dot{X}y ~\text{and}~ w'=y\dot{X} & 
\text{if}~d (X)=2\\
\text{when K means a cusp,} & w=Ky ~\text{and}~ w'=yK & \text{if}~d (K)=1\\
& w=\dot{K}y ~\text{and}~ w'=y\dot{K} & \text{if}~d (K)=2
\end{array}\right. \] 
This relation $\sim$ is called {\itshape marked cyclic equivalent for fronts}.  
\end{dfn}

Replacing cyclic equivalent for fronts by {\itshape marked cyclic equivalent for fronts}, we can easily  check that the similar results are established and can easily generalize Sect. \ref{f2-f3}.     
Therefore we only see the case of $FI_3$.  
Let $i$ be index.  
For every front $F$, 
we define $GFI_3$ by 

$
GFI_3 (F; x, y, z, p, q, r, s, t, u, v, h)
 := GCI_3 (F; x, y, z)
+\Big{[}p[XKXYY]+q[KXXYY]\\
+r[XKYXY] 
+s[KKK]
+t[XXKK] 
+u[\dot{K}XX]
+v[K\dot{X}\dot{X}]
+h[\dot{K}K], w_F \Big{]}.  $

\begin{cor}
$GFI_3$ is an invariant of fronts.  
\end{cor}

\begin{ex} 
There exist two fronts 
such that  the value of $FI_3$, index, Maslov index, basic invariants are the same; however, 
the value of $GFI_3$ on one curve is
 different from that on the other (Figure \ref{front4}).

\begin{figure}[htbp]
\begin{center}
\includegraphics[width=20cm,height=3cm,keepaspectratio]
{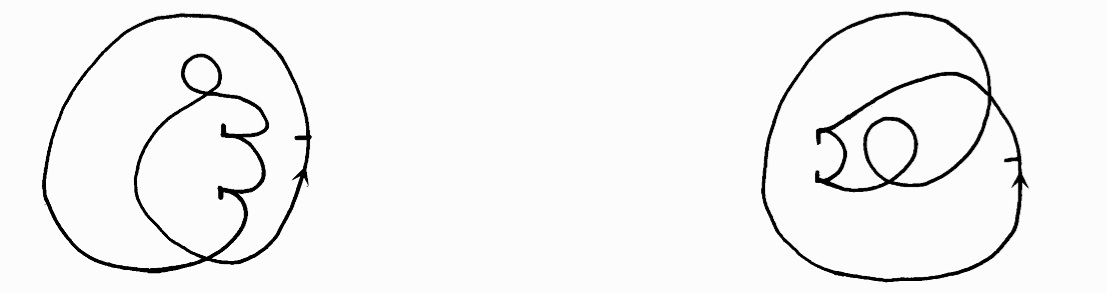}
\caption{Two fronts such that 
$($$FI_3,$ $i,$ $\mu,$ $J^+,$ $J^-,$ $St$$)$
 is $($$2q+2,$ $2,$ $0,$ $-4,$ $-5,$ $2$$)$. The value of  $GFI_3$ is 
$2 q - 2 u + 2 v + 2$ for the left figure,  $2 q + 2 u + 2$ for the right figure.  }
\label{front4}
\end{center}
\end{figure}
\end{ex}

In particular, this example implies the following.  

\begin{cor}
$GFI_3$ is a stronger invariant than  $FI_3.$  
\end{cor}

\section{Generalization}\label{strong}

By replacing 
 the function $\mathit{sign}$ with $\rho$, 
we get more general invariants.  

\subsection{Generalization of $I_n$}

Let $\tilde{\alpha} $ be commutative monoid generated by  an alphabet $\alpha$
 which contains the unit element and consider monoid algebra $\mathbf{Q}\tilde{\alpha}$.  
For every  $A \in \alpha 
$-alphabet $\mathcal{A}$
, suppose mapping  $\rho : |A| \to \rho (|A|) \in \mathbf{Q}\tilde{\alpha}$ is given.  

When an arbitrary \'{e}tale word  $w $ is given,  choose a sub-word  $v $ of  $w $, for this sub-word  $v $, 
sub-\'{e}taleword  $ (\mathcal{A}_v, v) $ of  $w $ is determined , and then  
\[\rho (v):= \prod_{A\in{A_v}} \{\rho (|A|)\}^{d (A)}\] can be defined.  
By using this, for every \'{e}tale word $w$
 and every sub-word $v$ of $w$, we define 
the mapping\[\left\langle v, w\right\rangle :=\sum_{\text{a sub-word $v'$ of $w$ isomorphic to $v$}} \rho (v').  \]  

Let $W_n$ be the $\mathbf{Q}$-module generated by the set  of 
all of the fake-Gauss words in $FW_n$.  
Let $FN_d$ be the free $\mathbf{Q}$-module generated by the set of fake nanowords over $\alpha$ $\{(\mathcal{A}, w)\}$ such that $card\mathcal{A}$ is
 less than $d+1.  $
Expanding
 $\left\langle v, w\right\rangle$  bilinearly, 
we can make a bilinear mapping$\left\langle, \right\rangle$ 
from $FW_n \times FN_d$
 to $\mathbf{Q}.  $    
For an arbitrary surface, 
let $w_\Gamma$ stand for a word which is determined
by a curve $\Gamma$ on the surface.  

\begin{thm}\label{I_n}
The following  $\{\tilde{I}_n\}$ (invariant of degree $n$) is the sequence of  surface isotopy invariants for
 pointed curves on a surface.  

\[\tilde{I}_n (\Gamma)=\left\langle \sum_k x_k v_k, w_\Gamma\right\rangle
\qquad  (n\in \mathbf{N})\]

where $\{v_k\}$ is the basis of $W_n$ and 
each $x_k $ is a parameter.  
\end{thm}

\begin{proof}
Theorem \ref{I_d} and the construction 
of $\langle , \rangle$ deduce immediately this theorem.  
\end{proof}

\subsection{Investigation of the generalization
 in the case $FI_2$}

For example, we will consider the case 
$FI_2$, and so 
we get an invariant $\widetilde{FI_2}$
 which is stronger than $FI_2$.  
For every element $A \in \alpha_{*}$-alphabet
 $\mathcal{A}$, 
$\rho$ is defined by $\rho (|A|)=|A| \in \overline{\mathbf{Q}\tilde{\alpha_{*}}}$.  

Let $F$ be front. By using
 the construction of $\left[, \right]$, 
the following  $\widetilde{FI_n}$ is a plane isotopy invariant of fronts.  
\[\widetilde{FI_n} (F)=\left[ \sum_{i,j} x_j [v_i], w_F \right]
\qquad  (n \in \mathbf{N}).  \]
$[v_i] $ in  
 $\overline{FW}_n $ (cf. Proposition \ref{[f]}),   
 each $x_j$ is a parameter.

For front $F= (\mathcal{A}_F, w_F) $, 
let $i$ be index, $2c$ the number of cusps, 
\[n_{\epsilon}=card\{A \in \mathcal{A}_F:|A|=a_{\epsilon}~
\text{or}~b_{\epsilon}\}
~(\epsilon = +,~\text{or}~-).  \]  

We define $\widetilde{FI_2}$ by 
\begin{eqnarray*}
\widetilde{FI_2} (F; p, q, x, z, t, v, r) &:=& pn_+ + qn_- + \langle
 x XXYY - x XYYX + z XYXY \\
& & + t KXX + t XXK - t XKX + vKK, w_F \rangle +rc \\
& & +\frac{x}{2}-\frac{x}{2}i^2.  
\end{eqnarray*}

\begin{thm}
$\widetilde{FI_2}$ is an invariant which is stronger than $FI_2$.  
\end{thm}

\begin{proof}

The value of $FI_2$ is obtain from $\widetilde{FI_2}$ by regarding $a_+$, $a_-$ as $1$.  
Therefore $\widetilde{FI_2}$ is at least
 as strong as $FI_2$.  
So,  
Figure \ref{front3} deduces that
 $\widetilde{FI_2}$ is an invariant which is stronger than $FI_2$.   

\begin{figure}[htbp]
\begin{center}
\includegraphics[width=14cm,height=3cm,keepaspectratio]
{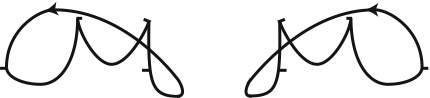}
\caption{The values of  $FI_2 $ for these two fronts are $p$ $+q$ $+x$ $+2t$ $-v$ $+2r$.  
The values of  $\widetilde{FI_2} $ for left figure is  $p + q + (-\frac{3}{2} + a_+ a_-) x + ({a_-}^2 + {a_+}^2) t + 
(a_+ a_-)v + 2 r $, 
The values of  $FI_2 $ for right figure is  $p + q + (-\frac{3}{2} + a_+a_-) x
 + ({a_+}^2 + 2 a_+a_- - {a_-}^2) t +
 (a_+ a_-) v + 2 r$.  }
\label{front3}
\end{center}
\end{figure}
\end{proof}

\begin{rem}
Suppose let  $F^{r} $ be the reflection of  $F $, and then 
 $\widetilde{FI_2} (F)=\widetilde{FI_2} (F^r) $.  
On the other hand, suppose
 $-F$ has only the difference of orientation from $F$, and then 
there exists an example (left figure in Figure 
\ref{front3}) as $\widetilde{FI_2} (F)\neq \widetilde{FI_2} (-F)$ (cf. $CI_2$, $LI_2$).  
\end{rem}

\begin{rem}
For these two fronts, $(i, \mu, J^+, J^-, St)= (2, 0,-2, -1, 2) $, and 
in terms of invariants of fronts :  $f^+, f^-, p^{\uparrow}, p^{\downarrow}, \lambda^{\uparrow}, \lambda^{\downarrow}
 $ due to Aicardi \cite{aicardi1}, 
$($$f^+,$ $f^-,$ $p^{\uparrow},$ $p^{\downarrow},$ $\lambda^{\uparrow},$ $\lambda^{\downarrow}$$)$ $=$ $($$2,$
 $0,$ $-2,$ $0,$ $2,$ $0$$)$.  
\end{rem}

\begin{rem}
 $\widetilde{FI_2} $ is the deformation 
$\overline{FI_2}$.  
Moreover, because  $\left\langle \dot{X}\dot{X}, w_F\right\rangle = n_ + {a_+}^2 + n_ - {a_-}^2$, if the term  $\left\langle \dot{X}\dot{X}, w_F \right\rangle $ is taken place of  $n_ + p + n_ - q $, the strength of the invariant does not change.  
\end{rem}

\vspace{0.3cm}
Department of Mathematics and Applied Mathematical Science Sciences

Graduate School of Fundamental Science and Engineering

Waseda University

3-4-1 Okubo, Shinjyuku-ku

Tokyo 169-8555, JAPAN

email: noboru@moegi.waseda.jp

\end{document}